\documentclass[11pt]{article}
\usepackage{amsmath}
\usepackage{url}
\usepackage{graphicx}
\usepackage{mwe}
\usepackage{amssymb}
\usepackage{hyperref}
\usepackage{amsthm}
\usepackage{authblk}
\theoremstyle{definition}
\newtheorem*{Hilbert90}{Hilbert Theorem 90}
\newtheorem*{Theorem}{Theorem}

\begin{document}
\title{A review of Pythagorean  Triples from both classical and  modern viewpoints}
\author{Ali Taghavi\\Qom University of Technology\\taghavi@qut.ac.ir}


\maketitle

\begin{abstract}
In this note we present a survey on some classical and modern approaches  on Pythagorean triples.   In particular some  approaches from non commutative  geometry,  operator theory and  also group theory  of Pythagorean triples are discussed.
\end{abstract}
\section*{Introduction}
A  Diophantine equation is  a  polynomial equation with integer coefficient for  which  only  integer solutions are  considered. The  Pythagorean equation \begin{equation}\label{Pythagoras}x^2+y^2=z^2 \end{equation} is one of the most famous  Diophantine equation. Every integer triple $(a,b,c)$ satisfying \eqref{Pythagoras} is called  a Pythagorean triple. In this paper we denote by PT as an abbreviation  for the set of Pythagorean triples. Based on the Pythagoras equation it is  constructed the Pythagoras tree initially presented by a  Dutch mathematics teacher  Albert E. Bosman in 1942. This tree has a fractal structure  whose s various fractal aspects including the spray  structure,  fractal  dimensions and  corresponding  zeta  and Laplacian operators can generates some lines of  researches. Regarding these  fractal terminologies   see \cite{FractalGeometry}.\\

The  Pythagorean theorem generates the concept of Landau Ramanujan constant: The  number of positive integers bounded by an integer $x$ which are the  sum of  two square is  asymptotically similar to $\frac{bx}{\sqrt{log(x)}}$ where $b$ is a constant  discovered by E. Landau and S. Ramanujan, see \cite{Landau}  and  \cite{Ram}

The Pythagoras equation  is  a  particular case of the equation $$x^n+y^n=z^n$$. Existence of   integer  solutions for this generalized equation is the Fermat's Last Theorem initially claimed in 1637 by Fermat and proved by Andrew Wiles in 1995 in \cite{Wiles}) states there are no positive integer solutions $(x, y, z)$. However the module $p$  version of this  generalized equation is studied in \cite{Shur} states that the Fermat last  theorem fails in finite  fields.  The Fermat last problem generates several  concepts in number theory. On of the most and  popular one is the Hardy  Ramanujan number 1729, the  least integer number can be represented in the form of sum of two cube number in two  different way: $1729=1^3+12^3=9^3+10^3$.\\
The PT have been known since ancient times. The oldest known record comes from Plimpton 322, a Babylonian clay tablet from about 1800 BC, in a sexagesimal number system. It was discovered by Edgar James Banks shortly after 1900, and sold to George Arthur Plimpton in 1922, for \$10. The name of  Pythagorean triples comes from the  classical  Pythagoras theorem stating that the  length sides $a,b,c$  of  a  right  triangle  with hypothenuse $c$ satisfy $a^2+b^2=c^2$ .  For  historical account on this  theorem  see \cite{STY}. Regarding the importance of the Pythagorean theorem and  also the time  credit of this theorem we quote from \cite{STY}:\\

\textit{"There are several outstanding contributions in mathematics that appeared in different cultures. The most notable is the Pythagorean theorem. Pythagoras is believed to be born between 580 to 568 BC. This is not much later than the time when the Egyptian papyrus was introduced into Greece around 650 BC. It was said that he  traveled to Egypt, Babylon and perhaps even India. Since there are ancient records of the Pythagorean triples in all these countries besides China, it is likely that Pythagoras learned the statement of his theorem from these countries. On the other hand, his great contribution is that he is the first one who proved the theorem. In fact, the concept of a proof based on formal logic found in Greek mathematics is rather unique as compared with other countries. It may be a mystery that countries far apart
discovered the same theorem at different levels of depth"}\\

 The  paper \cite{Givental}  contains a very  detailed  historical and  motivating  discussion. Historical considerations are also discussed  in \cite{Penrose}, \cite{RAVI}  and \cite{4000}. The latter explore a 4000 years history of the Pythagorean theorem.\\

   The paper  \cite{philosophy} has a philosophical discussion on PT. In this paper it is  mentioned that the  negative  conception of  mathematics developed  by  Hegel in his \textit{Phenomenology of  Sprit} is  based on the  study of  Euclidean proof  of the  Pythagoras theorem. Then the  author provide  serious critique on this  conception of  Hegel and express  that:\\

\textit{Not only are
there many other proofs, but what is at stake in Pythagoras’ theorem refers to
complex structures unnoticed by Hegel and which he could not know: relationship
between quadratic form and square of a linear form, geometric algebra, spinors and
rotations in the space, all concepts of great importance in modern physics. But these
are also very close, in fact, to what Hegel privileged: dialectical synthesis and movement. Thus, it is mathematics, and especially mathematical physics, which has now something hegelian, and maybe more than (contemporary) philosophy}\\

The  Pythagorean theorem is discussed very beautifully in \cite{Penrose} in particular in chapter 2 of the book. The chapter title is "An ancient  theorem and a modern question". The  chapter provides a geometric and illustrative proof  for the  theorem. The proof is based on covering of the plane with a family of  squares with two  different size. Then he provides another covering of the plane with a new family of  squares whose  vertices is the centroid of the bigger squares in the initial family. This consideration leads to an illustrative proof of the Pythagorean theorem.    Then it treat some well known proofs of the theorem and discuss about the implicit present of  parallel postulate in such proofs. Finally the chapter  presents some philosophical and  also physical investigation associated to this historical theorem  and  its  connection to other types of geometries.\\

This  amazing theorem in Euclidean geometry has a unification in two other types of geometries namely hyperbolic and spherical geometries. The  unification is the following theorem proved in \cite{Unification}:

\begin{Theorem} Let a  right triangle in a model of Euclidean, spherical and  hyperbolic geometry has sides $x,y,z$ where $z$ is the  hypotenuse. Then $A(z)=A(x)+A(y)-\frac{K}{2\pi}A(x)A(y)$ where $K$ is the curvature and $A(r)$ is the area of circle of radius $r$.

\end{Theorem}

This situation is a motivation to pose the question of integral PT in hyperbolic geometry.\\

 A  PT in particular represent three points on the plane whose mutual distances are integer. The generalization of this concept to higher number of points is presented in \cite{Anning}. In this paper P. Erdos and Norman Anning proved the following  result: An infinite number of points in the plane can have mutual integer distances only if all the points lie on a straight line. The  n dimensional version of this results is true as it is  indicated in the last paragraph of \cite{Anning}. This theorem generates the concept of Diophantin Erdos graph: A graph in the  plane which is  maximal with respect to the property that    all vertices   have integral  coordinates and  all distance are integers. Inspired by Erdos Anning theorem one may pose the following  question: Let we have a finite set of points $z_1,z_2,\ldots,z_n$ on the plan whose mutual  distances  are at least 1. Is there a straight line $L$  such that the set $\pi_L(z_1),\pi_L(z_2),\ldots,\pi_L(z_n)$  consist of  points  with integral mutual  distances? Here  $\pi:\mathbb{R}^2\to L$ is the  orthogonal projection. If the answer to this question is negative, inspired by the least square method in numerical analysis, one may pose the question of existence of a line which minimize   the deviation of the set of points $\pi_L(z_i)$ from being mutual integral distanced . More precisely assume that $A\subset \mathbb{R}^+$ is a finite set of positive real number. We associate to $A$ the  total deviation from integrality as $\sum_a (a-[a])^2$. Now one may pose the following question: For a  finite set $z_1,z_2,\ldots,z_n$ in the plane is there a line $L$ which minimize the total deviation of $A$ from integrality where $A=\{d(\pi_L(z_i),\pi_L(z_j)))\mid 1\leq i,j\leq n\}?$\\

If $(a, b, c)$ is a PT, then so is $(ka, kb, kc)$ for any positive integer $k$. A  primitive PT is  an integer  triple $(a,b,c)$ with $a^2+b^2=c^2$  and $gcd(a,b)=gcd(a,c)=gcd(b,c)=1$ namely they  are  mutually coprime integers. It is  well known that every primitive PT $(a,b,c)$  after a possible  exchange $a:=b, b:=a$ can be  written in the  form of  Euclid formula \begin{equation}\label{Pgenerator} a=r^2-s^2, b=2rs,  c=r^2+s^2\end{equation} for  coprime  integer $r,s$. For  a proof  see \cite{Mirzakhani}.\\

 The  similarity between this  representation and  component of  complex number $z^2=(x^2-y^2)+2xyi$  for $z=x+yi$ suggests that   at the time of Euclid, vague manifestations of the spirit of complex numbers appeared to mankind, but until it was fully revealed, mankind had to endure about 20 centuries.\\

We shall explain that the  Euclid formula  \eqref{Pgenerator} has  a miraculous relation to the  Hilbert theorem 90. To state the theorem we recall that  for  a finite field extension $K\subset L$ and  $a\in L$ the relative  norm $N(a)$ is  defined  as  the  determinant of  the  $K$-linear  map on $L$  with $x\mapsto ax$. For the particular case that the  Galois group is  a  finite  cyclic  group  generated by  $\sigma$ then $N(a)=a\sigma(a)\ldots\sigma^{n-1}(a)$  where $n$ is the order of $\sigma$. The  theorem says:

\begin{Hilbert90}
Let $K\subset L$ be  a Galois  field  extension. Assume  that the  Galois group $ Gal(L/K)$ is  a  cyclic group of  order $n$ generated by\newline $\sigma\in Gal(L/K)$. For  a  given  element $a\in L$  with relative norm $N(a)=1$ there exists $b\in L$  with $a=\frac{\sigma(b)}{b}$

\end{Hilbert90}

The theorem has a  generalization in terms of  Group  cohomology: For  every  Galois field extension $K\subset L$ the  first  group cohomology $H^1(G,L^*)$  vanishes. Here $G$ is  the  Galois  group $G(L/K)$  and $L^*$ is the  multiplicative  group of $L$.\\

Now  we explain that the  Hilbert theorem 90 implies the  Euclid representation of  PT. We apply the  Hilbert theorem 90 to  field extension $\mathbb{Q}\subset \mathbb{Q}(i)$, the  quadratic  extension of the field  of rational numbers. The  Galois  group is  a group of  order $2$  generated by conjugate automorphism $\sigma(x+yi)=x-yi$. The relative norm of  every  $a=x+yi$ is  $N(a)=a\sigma(a)=a\bar{a}=x^2+y^2$. So an  element  $a$ with $N(a)=1$ is  a rational  element $x+yi$ with $x^2+y^2=1$. Now the Hilbert  theorem  90 says that every rational point $x+yi$ on the  unit  circle  has a representation $x+yi=\frac{\sigma(b)}{b}$  for  some $b=c+di$  where $c,d$ are rational  number. After a possible  multiplication of  numerator and  denominator by a  natural number in  fraction $x+yi=\frac{\sigma(b)}{b}$ we  may assume that $c,d$ are  integers. So  we  have $$x+yi=\frac{c-di}{c+di}=\left(\frac{c^2-d^2}{c^2+d^2}\right)+\left(\frac{-2cd}{c^2+d^2}\right)i$$
This  coincides to Euclid representation \eqref{Pgenerator}.\\

The  situation mentioned  above is  a  motivation to pose the  following  question: Is it trivial to  compute some  higher  group  cohomologies $H^{i}(G, \mathbb{Q}(i)^*)$  where  $G$ is the  Galois  group of  the  field  extension? Are there  some  number theoretical interpretations for such a computation?\\

The formula \eqref{Pgenerator} ensure us that there  are infinitely many  solutions to this particular Diophantine  equation. The Euclid formula in particular implies that the area of every right triangle is an integer. This pose the question that what integers  arise as the area of a right triangle  with integer length sides. The question is explored in \cite{Ramin}. A related question to the latter is the following: Which integers can be the area of a right triangle with rational sides. It can be  shown that $1,2,3$ and $4$ can not be the area of an rational sided  right triangle. Determination of  all possible $N$ which can be the area of a rational sided right triangle is  closely related to existence of rational points on  elliptic curve $y^2=x^3-N^2 x$ see \cite{Tunneli} and \cite{Koblitz}.  The  triangle  with integer  length sides and integer  area  are called  Heronian  triangle. Obviously the  adjunction of two Pythagorean triangle along a common side  leads to a  Heronian triangle. The interesting  point is that the  converse is also  true. Namely every Heronian triangle is obtained in this way, see \cite{Carlson}\\
The right  triangle  with integer  sides and prime  hypotonus   is discussed in \cite{primeHypo}. In this paper they essentially use the  Hardy Ramanujan inequality to classify triangle with prime  hypotonus.
 There are several  algebraic  structures on the  space of  all primitive PT. For example the  multiplication  $$(x_1,y_1,z_1).(x_2,y_2,z_2)=(x_1x_2-y_1y_2,x_1y_2+x_2y_1,z_1z_2)$$ gives a Free abelian group structure  with infinite rank to the space of all primitive PT  module  similarity, see \cite{Eckert}.  A consequence  of this multiplication formula is the fame 2-square identity:$$(a^2+b^2)(c^2+d^2)=((ac-bd)^2+(ad+bc)^2$$ This  identity says that the  set of all positive integers who are the  sum of two square is closed  under integer  multiplication. The meaning of this identity is that the product of two  quadratic forms can be written in the form of composition of  two new quadratic forms. This  point of view and also similarity of the multiplication of PT to the complex  multiplication inspired by Gauss idea is  explained in \cite{linearalgebra}.\\
  A  base for the free Abelian group of PT is $$\{(a,b,p)|p^2=a^2+b^2, \quad \text{p=4k+1  is  a Pythagorean prime}\}$$ This is isomorphic to the space of all primitive triples module similarity of  triangle. The latter can be modeled as the group of rational points on unit circle, see \cite{Montly}. A monoid structure   is explained in \cite{Group}. Regarding the  group structure of PT
 it  would be interesting to produce some group theoretical invariant associated to this group  that at the same time these invariants would have a plane geometric nature. For examples one may pose the following questions:\\
  \textbf{Question   1}
  Is the  following quantity  an appropriate quantity as a group theoretical invariant? Is it invariant under all isomorphism of the group? To every Pythagorean triangle $\Delta$ whose edges forms a primitive triple  one associate the  maximum number of disjoint Pythagorean triangles can be embedded in $\Delta$.\\

  \textbf{Question  2} A Pythagorean triples is a particular case of an integer triangle: A triples of integers $(m,n,r)$ which are length of a triangles. One may define a primitive integer triangle analogously. So is there a  natural extension of  group structure of Primitive pythagorean triples to the  space of all primitive  integer triangles?\\

  A  ring structure  on the space  of  PT is  discussed in \cite{Missouri}. The ring  structure is defined  as  follows: Let we  have  two particular  PT $A=(a, \frac{a^2-n^2}{2n}, \frac{a^2+n^2}{2n})$  and $B=(b, \frac{b^2-n^2}{2n}, \frac{b^2+n^2}{2n})$ then the  operations  are
     \begin{equation*} A \oplus B=\left(a+b-n, \frac{(a+b-n)^2-n^2}{2n}, \frac{(a+b-n)^2+n^2}{2n}\right)\end{equation*}

      \begin{multline*} A \odot B=\bigg{(}(a-n)(b-n)+n, \frac{((a-n)(b-n)+n)^2-n^2}{2n},\\\frac{((a-n)(b-n) )^2+ n^2 }{2n} \bigg{)}\end{multline*}

      This  definition has  an extension to all  arbitrary  PT.\\

 The usual   strategy for  consideration of Diophantine problems is to investigate the problem of  existence of  rational points on corresponding algebraic  variety. As a notable case in this regard we point out to the  Norman Anning conjecture presented in early second decay of 20th centry. The  conjecture says that  it is possible to arrange $48$ points on a circle, such that all distances between the points are integer numbers and are all among the solutions of the diophantine  equation $x^2+xy+y^2=7^2.13^2.19^2.31^2$. A  geometric  approach to a generization of this conjecture is presented in \cite{aeqation}\\
    Apart  from Euclid original  formula mentioned  above  there  are  some  other formula  for  generating PT. One  interesting  method  is the Fibonacci  method. Consider  the  sequence of  odd positive  integers  $1,3,5,7,\ldots,$ The  sum of  $n$-th term of  this  sequence is  always  a perfect  square. Let $k$  be  an odd  perfect  square in the  above  list. Assume that  $k$ is  the  nth  odd  number in the list hence $n=\frac{k+1}{2}$. Put $a^2=k$. Let  $b^2$  be the sum of all $n-1$ first terms in the above list. Put  $c^2$ be the sum of  all nth terms. Then the  triple $(a,b,c)$ is  a primitive PT.  In the other words let $F_n$ be the  sequence of  Fibonacci numbers then the above  process generate a  PT $(F_nF_{n+3}, 2F_{n+1}F_{n+2},F_{n+1}^2)+F_{n+2}^2$. With this  construction the  consecutive  Fibonacci  numbers  are related via Farey sum of PT. Namely we  have  $F_n \oplus F_{n+1}=F_{n+2}$   where the operation $\oplus$ refers to the  Farey  sum of PT. For  Farey  sum see \cite{Farey}. Regarding  Farey  sum and  Farey diagrams, there is a beautiful geometric  description of  PT in terms of  Farey diagram in \cite{Hatcher}. In this  consideration the rational points on the unit circles corresponds to from one hand to PT from other hand to the rational real numbers via stereographic  projection. The Farey  diagram on the rational numbers is  constructed on the  nested  finite subsets of  rational numbers in [0,1] in the form $F_n$ consist of  all irreducible fractions $\frac{p}{q}$ with $p\leq n$. This  defines  a  graph  structure on whole  extended rational numbers $\hat{\mathbb{Q}}$. The  inverse  stereographic gives  us  a diagram structure on rational points on the  unit  circle. This situation relates to hyperbolic geometry since the unit  circles is the  boundary of the  disc model of hyperbolic geometry. The beautiful illustration is presented in \cite{Hatcher}. This  situation is  identical to the  Apolonian gasket mentioned in \cite{Clifford}\\
     As some other approaches   for  generating  PT we  point out to the Stifel approach: For  every mixed number $n+\frac{n}{2n+1}$  we represent  the corresponding  improper  fraction in the  form $\frac{b}{a}$ then the triple $(a,b,b+1)$ is  a primitive PT. For  example $2 \frac{2}{5}$ is simplified to $\frac{12}{5}$ so the  triple $(5,12,13)$ is a  primitive  PT. These approaches generate  some  of  such triples but not all of them. For  some  new  perspective  of  generating  of  all PT see \cite{NewPerspective}

 In \cite{philosophy}, apart from philosophical considerations,  there  are also  some    physical interpretation and  spinor representation  for these triples. For details on geometric and  physical aspects of  spinors see \cite{Bourguignon2}.  To have an informal but very inspiring introduction of  spinors we quote  from \cite{Bourguignon}the  following:\\

 \textit{In an interview with  Andrew Hodges in 2014 he said:"I was mystified by  spinors as they turn out to be square roots of  vectors and I  coud  not  understand how  you could do that"}.\\

  This  description would be  more  clear in the  following interpretation mentioned in \cite{philosophy}:  A  PT $(a,b,c)$ represents a  matrix in the  form $X=\begin{pmatrix} c+b&a\\a&c-b\end{pmatrix}$  whose  rank is  equal to $1$ hence there exist an integer   vector  $\psi=[p,q]^T$  with $X=2\psi \psi^T$  where $T$ is the transpose operation. Such a  vector  $\psi$ is  a  spinor. The  space of  spinors  can be  acted  by the  modular  groups. The  modular  group is the  quotient of  the  group $Sl(2,\mathbb{Z})$ by $\pm I_2$  where $Sl(2,\mathbb{Z})$  is the  group of  all 2 by 2 matrices with integer  coefficient with determinant equal to 1. The  action of  modular group on the space of  spinors  gives  an  action on the  space  of  PT. The  action is defined  via  $A.X=AXA^T$. The action  fails to keep invariant the  space of  primitive PT. This  leads  the  author of \cite{philosophy} to define the standard primitive PT as  follows: A triple $(a,b,c)$ is called a  standard triple if $c>0$  and either $(a, b, c)$ are relatively prime or $(a/2, b/2, c/2)$ are relatively prime with $a/2$ odd. With this  concept of  standard  PT one  observe that the  space  of  standard  triples  is invariant under the  action of  modular  group.  The  computation in \cite{philosophy}  turn out that the  space of  primitive  PT is  invariant under the special subgroup of  the  modular  group obtained as the  kernel of the  natural group transformation $Sl(2,\mathbb{Z})\to Sl(2,\mathbb{Z}_2)$ where $\mathbb{Z}_2$ is the  field of order $2$ namely the  integers module $2$. This  implies that every  primitive triples is  uniquely representable in the  form of  products of two  matrices $U=\begin{pmatrix} 1&2\\0&1 \end{pmatrix}$ and $L=\begin{pmatrix} 1&0\\2&1 \end{pmatrix}$ \\

We  note that the  matrix  $U$ above is the  standard  example of  every  books  on hyperbolic  and  Anasov  toral  diffeomorphism.  So  it  would be interesting to  find  a possible  dynamical relation between PT and torus Anasov and  hyperbolic diffeomorphisms.See \cite{Anasov}\\

In \cite{Clifford}  a  more detailed  spinor representation of PT is explained. In the  latter it is deeply introduced  several  geometric  approach to  Pythagorean triples via  symplectic  geometry, clifford  algebra  and  also  Apollonian gaskets. In particular it asserts that every  primitive PT is  represented in the  Apollonian window.\\

The  dynamic of   PT represented as rational points on the unit circle is  presented in \cite{Dynamic}. In this paper certain ternary expansion is associated to a rational point of  unit circle then based on this expansion  some ergodic and dynamical information is extracted. This  ternary expansion is closely related to the  rooted tree  structure of PPT with a unique   generating  node $(3,4,5)$  which is attributed to  B. Berggren \cite{Bergern}  and J.M. Barning,  \cite{Barning}.  In fact the PT are  integer solution of  certain integral  quadratic  form $x^2+y^2-z^2=0$ and the latter two references provid  tree structures for the  solutions of this  particular integral ternary  quadratic form. This  situation is  generalized to other types integral  quadratic  forms in \cite{Journalofnumbertheory}.\\
As  another  dynamical  aspects of such triples  we  point out to the  the  billiard structure on primitive PT. This  structure  is studied in \cite{Panti}. In this  consideration the  billiard maps act on the  space $x_1^2+x_2^2-x_3^2=1$. The restriction to the unit circle has a factor $B$. Then the primitive PT are $B$ preimage of the fixed points of the billiard.

There is  a  beautiful  description and  generalization of PT in terms of  Riemannian  Geometry. Lets generalize the  equation $x^2+y^2=z^2$ to the  matrix  form $A^2+B^2=C^2$  where $A,B,C$ are matrix of  order $n$.  Let $M^3(c)$ be  the  standard  complete  simply  connected Rimannian manifold with  constant  sectional curvature $c$  where  $c\in \{-1,0,1\}$. Assume that $M_2$ is  a two  dimensional manifold  which is  isometrically  immersed in  $M^3(c)$. So it naturally generates the  first, the  second  and the  third  fundamental  forms $\textbf{I}, \textbf{II}$ and  $\textbf{III}$ We say  that $M_2$ is a Pythagorean  surface if  these 3  matrices  satisfy the  Pythagorean relation $\textbf{I}^2+\textbf{II}^2=\textbf{III}^3$. This  concept  can be  generalized  to  arbitrary  dimension. In this  setting it is proved in \cite{Riemann}  that every Pythagorean  Hypersurface in $\mathbb{R}^{n+1}$ is  isometric to the  round  sphere of  radius $\phi$ where $\phi$ is  the  conjugate  golden number $\phi=\frac{-1+\sqrt{5}}{2}$. In \cite{bounded} a comparison version of Pythagorean theorem  is presented  to judge the lower or upper bound of the curvature of Alexandrov spaces\\ Pythagorean  submanifolds  isometrically immersed in Euclidean spaces is  studies in  \cite{Pythagorean triple inEucleadian space} \\

A  complex  algebro geometric  approach to the  Pythagorean theorem is presented in \cite{complex triangle}. The author present the  concept of  complex triangle, a triangle whose  sides can be  complex number and present an appropriate Pythagorean theorem.

There is  a non commutative analogy of Pythagoras theorem mentioned in \cite{NCG}. The  usual Pythagoras theorem can be  formulated in terms of Connes distance on the  space of  pure  states in the  product of  commutative  spectral triples. This commutative  formulation is generalized  in \cite{NCG} to non commutative  spectral triples and (not  necessarily pure) states. In this way  a Pythagorean inequality is discussed. With this  non commutative  view point on Pythagoras theorem one may think to an  appropriate  analogy of PT namely the 3  distance between the underling (pure) states would be integers. Is there a possible relevant relation between this question and  the idea of  natural Dirac operator mentioned in \cite{Dirac} ?\\

A very inspiring  number theoretical  approach to  right triangles  and  PT is presented in \cite{Ramin}. The  book  explores the  distribution of  rational points  on unit  circles. It also provide an alternative  proof  of the  four  square theorem.\\

The  PT has been studied also from a  Ramsey theoretical point of view. The Boolean PT problem asks  if it is possible to color each of the positive integers either red or blue, so that no PT $a,b,c$  have the same  color. Such a  coloring is possible up to $7824$. But for every coloring of  $1,2,\ldots,7825$ with two colors there is a PT in this interval all have the same color, see \cite{color}.\\

The  distribution of  PT  are  discussed in \cite{distribution}: For  every integer $x$  we  denote by $D(x)$ the  number of  all $(a,b,n)$  of  triples  with  hypotenuse $n$  such that $n\leq x$. The  asymptotic  expansion of this  function is  discussed in \cite{distribution}. The  paper introduced a relation between this  asymptotic  expansion and distribution of  zeros of  the Riemann  zeta  function.

\textbf{ Remark   1}. The Pythagorean theorem states that $x^2+y^2-z^2=0$  for  leg lengths $(x,y,z)$ of a right triangle. It generates in particular the  problem of existence and  distribution of Pythagorean integers. Note that $x^2+y^2-z^2=0$ defines  an  algebraic curve. Is it  a  good  idea to  peak another  quadratic integral  form say $kxy=z^2$ then pose the analogue questions in this new setting? It can be  shown that a triangle with sides length $(x,y,z)$ whose  side  with length $x,y$ make an angle $\theta$ satisfies this  algebraic  equation if  and only if $cos \theta= \frac{k}{2}\times \frac{x^2+y^2-z^2}{z^2}$. So  a  right  triangle case is a very particular case of this  situation. What is the advantage of  $x^2+y^2-z^2=0$ compared to the new  algebraic curve $kxy=z^2$?  As two  complex  algebraic curve they  are  equivalents. But what  about their  number theoretical  aspects?\\

\textbf{Remark   2}  There is  a non spliting  short exact  sequence
  \begin{equation}\label{SP}o\mapsto \mathbb{Z}^3\to P\to\mathbb{Z}_2\times \mathbb{Z}_2 \to 0\end{equation} for  which the  group $P$ plays a crucial role in the  unite  version of  the  Kaplanski  conjecture. For  such group $P$  the  group  algebra $\mathbb{C}P$  has  non trivial unit. This  is  a  counter example to the unit  conjecture of  Kaplanski provided by  Gill  Gardam, see \cite{Gill} Recall that the  Kaplansky unit conjecture says that for a torsion free group $G$ and a field $F$ the  group ring $FG$ has no nontrivial unit. Gill Gardam first presented a  counter example for $F=\mathbb{Z}/2\mathbb{Z}$ then  generalized to the  field of  complex numbers. So it would be interesting to find a  number theoretical interpretation for this  group $P$. The  above  short  exact  sequence suggest an embedding of  $\mathbb{Z}^3$ in group $P$. In particular the  group $P$ contains all PT. The  first  question is that does the  set of  all PT generate the  group $\mathbb{Z}^3$? Is it a  base for this  group?If the answer to the  both question is  affirmative then one may define  an area  function on $P$ as  follows: The  group $P$ is the disjoint  union of $\mathbb{Z}^3$ and its two cosets according to splitting in \eqref{SP}. Since $\mathbb{Z}^3$ can be  generated by all PT then one may define the area functions  via linearity then extend to the other two cosets of $\mathbb{Z}^3$. This  area function is not  an element of  $\mathbb{C}P$ since it is not  a  finite  support  function. But is it a point wise limit  of  elements  $u_n$ of  $\mathbb{C}P$  such that every $u_n$ is  a  unit  element of the group  algebra?\\

  \textbf{Remark  3} A pair $(m,n)$ of integers is called  a Pythagorean pair if $\sqrt{m^2+n^2}$ is an integer. The Euclid representation of PT can be read as follows: The quadratic function $f:\mathbb{C}\to \mathbb{C}$  with $f(z)=z^2$ maps the the integer Latic $\{m+ni\mid m,n\in \mathbb{Z}\}$ to the  Pythagorean pairs. So for every integer $k$ the monomial $kz^{2^n}$ has this property. So one may pose the  following  question:  are there other type of polynomials in one complex variables or other type of entire functions which maps the Lattice of integers to the set of Pythagorean pairs?\\
  
  \begin{figure}
\centering
 \includegraphics{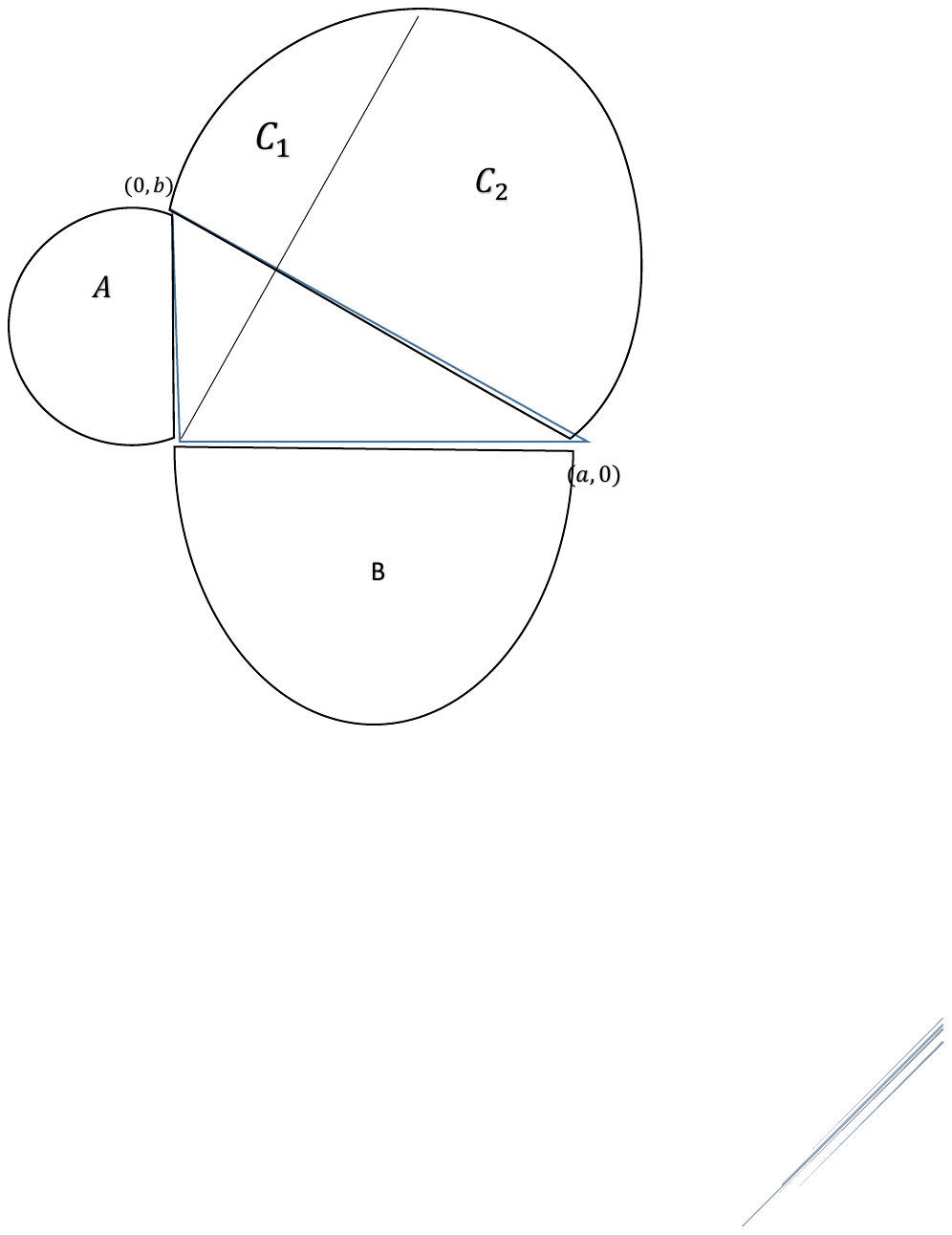}
 \end{figure}

  \textbf{Remark  4} One may pose the question of PT in a  $C^*$  algebra. An integer element of a  $C^*$ algebra $A$ is a self adjoint element $a\in A$ such that $\phi (a)\in \mathbb{Z}$ for all pure states $\phi$ on $A$. A PT is a triples $(a,b,c)$ of not necessarily commuting integers of $A$ satisfying $a^2+b^2=c^2$. A scalar triple is one which all components are scalar. So under what conditions does a $C^*$ algebra  not admit non scalar PT?

  \textbf{Remark   5} The classical  Euclid proof of the Pythagorean theorem is based on the following: he  considered a right triangle with length sides $(a,b,c)$   with $c$ hypotenuse. He built three squares $(A,B,C)$ on these three sides. He drew the height perpendicular to the hypotenuse. This height divide the square $C$ into rectangles $C_1,C_2$. Then he proved that the area of  $C_1$ is equal to the area of $A$ and the area of $C_2$ equals to the area of  $B$. This results to $c^2=a^2+b^2$. But the later can read in the equivalent form  $\pi c^2=\pi a^2+\pi b^2$. This is  a motivation to think to an alternative proof  for the Pythagorean theorem: One build three semi circle $A,B.C$on the  edges $a,b,c$ whose  diameters are $a,b,c$ respectively.
  Then one draw the  height to the  hypotenuse.Then pose the question does the height divide the  semi circle $C$ in to tworegions $C_1,C_2$ with the property that the area of  $C_1$ is equal to the area of  $A$ and the area of  $C_2$ is equal to the area of $B$?This would gives us an alternative proof  for the Pythagorean theorem. But in reality the answer to the  question is not affirmative in general unless the triangle satisfies $a=b$. So a  natural question is the following: Let we  have  a right triangle with sides $a,b,c$  with corresponding  vertexes $V_a,V_b,V_c$. Then there is a unique point $P$ on the  hypotenuse $c$  with the property that the segment $\bar{V_cP}$ divides the  semi  circle $C$ into regions $C_1,C_2$  with area of  $C_1$ is equal to the area of $A$ and the area of $C_2$ equals to the area of $B$. What is the  precise formulation of this  unique point in terms of vertex  coordinates? What is the properties of this point $P$ in terms of terminologies of classical and  ancient plan geometry?
  The  mathematica codes for our \textit{Mathematica} computations are the  following. Here we  consider a  right triangle with vertex $(0,0),(a,0),(0,b)$
  In our computations $S-w$ is the area of $C_1$ and  $c$ is the area of  semi circle $A$\\
  
  $t = VectorAngle[{-a/2,
  b/2}, {(-a/2) + 2 a*b^2/(a^2 + b^2), (-b/2) +
     2 a^2*b/(a^2 + b^2)}];\\
S = (t/8)*(a^2 + b^2);\\
w = Area[Triangle[{{a/2, b/2}, {2 a*b^2/(a^2 + b^2),
      2 a^2*b/(a^2 + b^2)}, {a*b^2/(a^2 + b^2), a^2*b/(a^2 + b^2)}}]];\\
c = Pi*b^2/8;$



\end{document}